\documentclass{aptpub}
\usepackage[square,compress,comma, numbers,sort]{natbib}
\usepackage[colorlinks=true, citecolor=red, linkcolor=blue]{hyperref}
\usepackage{amsfonts,mathtools}
\usepackage{amsfonts,mathtools,natbib}
\usepackage[shortlabels]{enumitem}

\allowdisplaybreaks[4]

\def\cadlag{c\`adl\`ag}

\usepackage{amsmath}
\usepackage{bbm}
\usepackage{amssymb}
\usepackage{mathtools}

\numberwithin{equation}{section}

\usepackage{color}

\definecolor{c20}{rgb}{0.,0.7,0.}
\definecolor{c30}{rgb}{0.,0.,1.}
\definecolor{c40}{rgb}{1,0.1,0.7}
\definecolor{c50}{rgb}{1,0,0}
\definecolor{c60}{rgb}{1,0.9,0.1}
\definecolor{c70}{rgb}{0.50,1.00,0.00}

\allowdisplaybreaks[4]

\usepackage{amssymb}
\usepackage{color}
\newcommand{\abs}[1]{\left\lvert #1 \right\rvert}

\DeclarePairedDelimiterXPP\pk[1]{\mathbb{P}}\{ \}{}{ #1}
\DeclarePairedDelimiterXPP\Ex[1]{\mathbb{E}}\{ \}{}{	#1}

\def\FRE{\mbox{Fr\'{e}chet }}

 \usepackage{xparse}% http://ctan.org/pkg/xparse

\NewDocumentCommand{\ceil}{s O{} m}{%
  \IfBooleanTF{#1} % starred
    {\left\lceil#3\right\rceil} % \ceil*[..]{..}
    {#2\lceil#3#2\rceil} % \ceil[..]{..}
}
\NewDocumentCommand{\floor}{s O{} m}{%
  \IfBooleanTF{#1} % starred
    {\left\lfloor#3\right\rfloor}
    {#2\lfloor#3#2\rfloor}
}

\newcommand{\norm}[1]{\lVert #1 \rVert}

\definecolor{c20}{rgb}{1.,1,0}
\definecolor{c30}{rgb}{1.,1,0}

\def\bE#1{#1}

\def\kk#1{{\textcolor{cyan}{#1}}}
\def\kkk#1{{\textcolor{cyan}{#1}}}
\def\kkkk#1{{\textcolor{cyan}{#1}}}

\def\kk#1{#1}
\def\kkk#1{#1}
\def\kkkk#1{#1}
\def\k2#1{#1}
\def\kdd#1{#1}
\def\ve{varepsilon}

\def\cL#1{{\textcolor{c40}{#1}}}
\def\cL#1{#1}

\def\vEE#1{{#1}}
\def\cEE#1{{\textcolor{c40}{#1}}}
\def\cEE#1{#1}
\def\rrd#1{\textcolor{red}{#1}}
\def\rrd#1{#1}

\def\tE#1{{\textcolor{c40}{#1}}}

\def\tE#1{#1}

\newcommand{\prooftheo}[1]{ \textsc{Proof of Theorem} \ref{#1} }
\newcommand{\proofprop}[1]{\textsc{Proof of Proposition} \ref{#1}}

\newcommand{\proofkorr}[1]{\textsc{Proof of Corollary} \ref{#1}}

\newcommand{\QED}{\hfill $\Box$}

\newcommand{\COM}[1]{}

\def\IF{\infty}

\newcommand{\R}{\mathbb{R}}
\newcommand{\inr}{\in \R}

\newcommand{\BQN}{\begin{eqnarray}}
\newcommand{\EQN}{\end{eqnarray}}
\newcommand{\BQNY}{\begin{eqnarray*}}
\newcommand{\EQNY}{\end{eqnarray*}}

\newcommand{\limit}[1]{\lim_{#1 \to   \infty}}

\newtheorem{theo}{Theorem}[section]
\newtheorem{korr}{Corollary}[section]

\def\bqny#1{\begin{eqnarray*} #1 \end{eqnarray*}}
\def\bqn#1{\begin{eqnarray} #1 \end{eqnarray}}

\newcommand{\BS}{\begin{sat}}
\newcommand{\ES}{\end{sat}}
\newcommand{\BT}{\begin{theo}}
\newcommand{\ET}{\end{theo}}
\newcommand{\BKK}{\begin{korr}}
\newcommand{\EKK}{\end{korr}}

\newcommand{\BEX}{\begin{example}}
\newcommand{\EEX}{\end{example}}

\newcommand{\BD}{\begin{de}}
\newcommand{\ED}{\end{de}}

\newcommand{\BIT}{\begin{itemize}}
\newcommand{\EIT}{\end{itemize}}
\newcommand{\BDI}{\begin{description}}
\newcommand{\EDI}{\end{description}}

\newcommand{\BRM}{\begin{remark}}
\newcommand{\ERM}{\end{remark}}

\newcommand{\BEL}{\begin{lem}}
\newcommand{\EEL}{\end{lem}}

\newcommand{\neprop}[1]{{Proposition \ref{#1}}}
\newcommand{\netheo}[1]{{Theorem \ref{#1}}}

\newcommand{\nekorr}[1]{{Corollary \ref{#1}}}

\def\TT{\mathcal{T}}
\def\KK{K}

\def\Z{\mathbb{Z}}
\def\inn{\in \mathbb{N}}

\newcommand\ind[1]{\mathbb{I}{\left\{#1\right\}}}

\def\TTT{\mathcal{T}}
\authornames{K. D\c{e}bicki, E. Hashorva, Z. Michna} % insert the authors here for use in running head
\shorttitle{On the continuity of Pickands constants} % insert short title here for use in running head

\begin{document}

\title{On the continuity of Pickands constants}

\authorone[University of Wroc\l aw]{Krzysztof D\c{e}bicki} % Affiliation is just the name of your university or institution

\addressone{Mathematical Institute, University of Wroc\l aw, pl. Grunwaldzki 2/4, 50-384 Wroc\l aw, Poland.

\email{Krzysztof.Debicki@math.uni.wroc.pl}} % Your postal address goes here.

\authortwo[University of Lausanne]{Enkelejd  Hashorva} % Affiliation is just the name of your university or institution

\addresstwo{Department of Actuarial Science,
University of Lausanne, UNIL-Dorigny, 1015 Lausanne, Switzerland.

\email{Enkelejd.Hashorva@unil.ch}}

\authorthree[Wroc\l aw University of Economics and Business]{Zbigniew Michna}

\addressthree{Department of Logistics, Wroc\l aw University of Economics and Business,
Komandorska 118/120, 53-345 Wroc\l aw, Poland.

\email{zbigniew.michna@ue.wroc.pl}}

\begin{abstract}
\rrd{
For a  non-negative separable  random field
  $Z(t), t\in \R^d$ satisfying some mild  assumptions we show that
\bqny{ H_Z^\delta =%=  \limit{n} \frac{1}{n \delta } \E*{\sup_{1 \le k  \le n} e^{ X(k\delta)}  }
		\limit{T} \frac{1}{T^d} \Ex*{\sup_{ t\in [0,T]^d \cap \delta \Z^d } Z(t) } \tE{<\infty}
}
for $\delta \ge 0$ where $0 \Z^d\coloneqq\R^d$
\kkk{and} prove that $H_Z^0$ can be approximated by $H_Z^\delta$ if  $\delta$ tends to 0.
These results
\kkk{extend the classical findings for the Pickands constants
$H_{Z}^\delta$, defined for
$Z(t)= \exp\left( \sqrt{ 2} B_\alpha (t)- \abs{t}^{2\alpha  }\right), t\inr$ with $B_\alpha$ a
standard fractional Brownian motion with Hurst parameter
$\alpha \in (0,1]$.} }
 The continuity of $H_{Z}^\delta$ at $\delta=0$ is additionally shown for two particular extensions of Pickands constants.
\end{abstract}



\section{Introduction}
 The  discrete Pickands constant $H_{Z}^\delta$ is defined for a given positive $\delta$ by
 \bqny{ H_{Z}^\delta=  \limit{n} \frac{1}{n \delta } \Ex*{\sup_{1 \le k  \le n} Z(k\delta)  }
 =
 	\limit{T} \frac{1}{ T } \Ex*{\sup_{ t\in [0,T] \cap \delta \Z }  Z(t)  } \in (0,\IF),
 }
 where $Z(t)= \exp(\sqrt{ 2} B_\alpha(t)- \abs{t}^{2\alpha})$ and $B_\alpha$ a standard
 fractional Brownian motion (fBm) with Hurst parameter $\alpha\in (0,1]$.
  When $\delta=0$, interpreting $0\Z$ as $\R$, $H_{Z}^0$ can be defined \kk{in an analogous way,} i.e., % and we have
 $$H_{Z}^0= 	\limit{T} \frac{1}{ T } \Ex*{\sup_{ t\in [0,T] } Z(t)  } \in (0,\IF),$$
 which is the classical  Pickands constant appearing in the tail
 \kkk{asymptotics of  the  distribution of supremum for a wide class of} Gaussian  processes,
 see e.g.,  \cite{PickandsB,Pit96,DiekerY}.
 Pickands pioneering method (see \cite{PickandsB}) for the approximation of \kkk{the tail distribution of supremum  for}
 a stationary Gaussian process relies strongly on a discretisation approach.
 	A crucial element in Pickands methodology is the fact that
 \bqn{\label{jojo}
 	\lim_{\delta \downarrow 0} H_{Z}^{ \delta}= H_{Z}^0.
 }
Notably, the first attempt of Pickands to prove  \eqref{jojo} contains  a gap; a correct proof is given in
\cite{MR0348906}[Thm B3], see also the comments after \cite{MR0348906}[Lem A3].\\
A systematic study of Pickands constants started with \cite{DiekerY}[Thm 1, Prop 2] showing that
\bqn{ \label{DMY}
	 H_{Z}^0 = \Ex*{ \frac{ \sup_{t\in \R} \kk{ Z(t)} }{  \eta \sum_{t\in \eta \Z } Z(t) }} = \Ex*{ \frac{ \sup_{t\in \R} \kk{ Z(t)} }{  \int_{t\in \R } Z(t) \lambda(dt)}}
}
is 	valid for any $\eta > 0$, where
$\lambda(\cdot)$ is the Lebesgue measure on $\R$. \\
%\kkk{In fact, the second formula in \eqref{DMY} is a consequence of	an earlier result obtained in \cite{Genna04}[Thm 2.1], see Remark \ref{remGenna} below.}\\
Both expressions in \eqref{DMY}  paved the way for simulation of Pickands
constants and
\kkk{inspired}  new  formulas for extremal indices of stationary time series,
see e.g., \cite{DM, Htilt, SBK, KDEH1, Hrovje, HBernulli,PH2020}.% where \kkk{related} expressions  are derived.

Pickands constants  for general Gaussian processes were considered first in \cite{Deb99}, see also \cite{MR2222683}.
Later \cite{debicki2017approximation, Htilt} discussed  extensions to general random fields (rf's) as we briefly outline next.
Let therefore $Z(t),t\in \R^d$ be a separable, non-negative rf  such that
\bqn{\Ex{ Z(t) }=1,  \quad \forall t\in \R^d
\label{MS}
}
and define  the corresponding Pickands constant  by
\def\HZA{H_{Z }^\delta  }
\bqn{\label{pcf}
\HZA  = 	\limit{T}T^{-d} \Ex*{\sup_{ t\in [0,T]^d \cap \delta \Z^d } Z  (t)  }, \quad \delta\ge0,
}
where we set $\delta \Z^d\coloneqq\R^d$ when $\delta=0$. If $\delta=0$, in order to avoid degenerated cases, we shall  assume further that
\bqn{
	 \label{norm}
\Ex*{ \sup_{t\in[0,T]^d} Z (t) } \in (0, \IF), \quad \forall\, T>0.
}
Definition (\ref{pcf}) might not be valid for a general $Z$ \cEE{since the limit  might not exist.}
\kkk{However,} if for all  compact sets   $K \subset \R^d$
\bqn{ \Ex*{ \sup_{t \in K }  Z (t+c)}= \Ex*{ \sup_{t \in K}  Z (t)} , \quad \forall c\inr^d,
	\label{shi}}
then as we shall show in Section 2   the constant  $\HZA  $ is well-defined and finite for any $\delta \cEE{\ge} 0$.
\kkk{Notably,} we know from \cite{Htilt,debicki2017approximation, HBernulli} that the limit in (\ref{pcf})
exists if  there is  a  stationary max-stable %\cEE{stochastically continuous}
process  $Y(t),t\inr^d$ with unit \FRE
marginals which has spectral process $Z$ \rrd{in} its de Haan representation  (\kkk{see} e.g., \cite{deHaan,dom2016})
\bqn{\label{eq1}
	Y(t)=  \max_{i\ge 1} \Gamma_i^{-1}  Z^{(i)}(t), \quad t\in \R^d.
}
Here $\Gamma_i= \sum_{k=1}^i \mathcal{V}_k$ with $\mathcal{V}_k, k\ge 1$ mutually independent
unit exponential rv's
being independent of $\{Z^{(i)}\}_{i=1}^\infty$, which are independent copies of $Z$.
The finite dimensional distributions (fidi's) of $Y$ are given by
\bqn{\label{frb}
\pk{Y(t_1)\leq x_1,\ldots, Y(t_n)\leq x_n} %=H_{t_1 \ldot t_n}(x_1,\dots, x_n)
=e^{- \Ex{ \max_{i\le n } Z (t_i) /x_i}}, \quad x_i>0, t_i \in \R^d, i\le n
}
and hence
$$ \left(\pk{Y(t_1)\leq x_1,\ldots, Y(t_n)\leq  x_n} \right)^m
=  \pk{mY(t_1)\leq  x_1,\ldots, mY(t_n)\leq  x_n}  $$
for all $m>0$,
which shows that the fidi's of $Y$ are  max-stable.
\cEE{Clearly, if $Y$ is stationary, then  $\sup_{t \in K }  Y (t+c)$ has the same law as   $\sup_{t \in K }  Y (t)$ for all $c\inr^d$ and hence if $Y$ has locally bounded sample paths,  then \eqref{frb} implies \eqref{norm} and further  \eqref{shi} is valid. }
\\
Pickands constants are \kk{closely} related to extremal indices  of the max-stable stationary  rf  $Y$.
\kk{Indeed,}   under the assumption \eqref{norm} and the finiteness of Pickands constants,
the separability of $Y$ implies
\bqn{\label{LY}
	-   \ln  \pk*{ \sup_{t \in [0,T]^d \cap  \delta \Z^d } Y(t) \le  r T^{d} } =  \frac{1}{rT^d}\Ex*{  \sup_{t \in [0,T]^d \cap  \delta \Z^d }  Z (t) }    \to  \frac 1 r \HZA
}
as $T\to \IF$ for all $r>0$. \kkk{Thus,} by definition,
the extremal index of the stationary rf $Y(t), t\in \delta \Z^d$ is equal to $ \delta^d  \HZA \in [0,1]$ for any $\delta>0$.
Clearly, \eqref{LY} is an approximation of  the
\kk{distribution of} \cEE{supremum of  $Y$. Such approximations are known for
general \cEE{stationary  rf's}. A prominent instance is $Y$ being a
symmetric $\alpha$-stable  rf, see \cite{Genna04,Genna04c,Roy,Roy0} and the references therein.
}
\\
\kkk{A natural question that arises here is the relevance of general Pickands constants, both in extreme value theory of
rf's and in stochastic modelling.}
As shown in \cite{MR3877258},  Pickands  constants  determined by
\kkkk{
$Z(t)= \exp( W(t)- \sigma_W^2(t)/2)$,
where $W$ is a centered
Gaussian  rf
with stationary increments and variance function $\sigma_W^2$,}
appear naturally in risk and queueing theory. Moreover,  as \kkk{advocated} in \cite{HBernulli},
Pickands constants related to a general non-Gaussian rf $Z$ have appeared in the literature in numerous
papers. \kkk{In that context,} considering a general  rf  $Z$ is important
since it unifies the study of extremal indices and Pickands constants.\\
\kkk{The second question is: for what  general $Z$ does the limit \eqref{jojo} hold? The answer to this question is presented in Section 2.}
Such a result has two immediate consequences, namely:
\begin{enumerate}[A)]
	\item  the discretisation method of Pickands for the study of  extremes of Gaussian rf's can be utilised also for cases where the limiting constants are determined by general Gaussian rf's with stationary increments;
\item  the calculation of   $H_{Z}^0$ can be carried out  by simulating  $H_{Z}^\delta$ for small $\delta>0$. \kkk{Therefore, it is further interesting to derive tractable formulas  for  $H_{Z}^\delta$ as given by \eqref{drejt}.}
\end{enumerate}

Organisation of the rest of the paper is the following.
In Section 2 we \kk{first} discuss the finiteness and the continuity of discrete Pickands constants  (in \netheo{th1})
and then   \rrd{\kkk{derive a} formula corresponding to \eqref{DMY} (in Proposition \ref{TM2}),}
which shows in particular how to approximate $H_{Z }^0$ using the fidi's of $Z$.
Section 3 is concerned with two extensions. The first one is motivated by results related to symmetric $\alpha$-stable
rf's \rrd{derived in} \cite{Genna04,Genna04c}, whereas the second extension is motivated by the following  constant
\bqn{ \label{hexh}
	\mathcal{H}_{Z}^\delta =   \limit{T} T^{-d} \int_0^1 \Ex*{\sup_{ t\in [0,T]^d \cap \delta \Z^d } Z_z(t)   }dz
}
defined in \cite{MR2206343} for  $\delta =0$. Here $ Z_z, z\in [0,1]$ is equal in law to $ \exp( \rrd{W}_z(t) - \sigma_z^2(t)/2),$ with \rrd{$W_z$}  a centered Gaussian  rf with stationary increments, variance function $\sigma^2_z$ and almost surely continuous sample paths.
 In \netheo{th11} we establish  the continuity of $\mathcal{H}_Z^\delta$ at $\delta=0$. Such a result is crucial for the applications of Pickands discretisation method  in the study of extremes of locally stationary Gaussian  rf's.  It is also of certain relevance  for the simulations of those constants. We have relegated all the proofs to Section 4.

\def\HDM{ \mathcal{H}_{X,r}^\delta}
\def\HDO{ \mathcal{H}_{X,r}}

\section{Main Results}
Our first result establishes the finiteness of Pickands constants and \eqref{jojo} for a general  rf $Z$ under some weak restrictions. In particular our result is satisfied for $Z$ such that $Y(t),t\inr^d$ defined in \eqref{eq1} is a stochastically continuous stationary max-stable  rf with locally bounded sample paths. \cEE{Hereafter  we shall suppose that all rf's are defined on a complete probability space $(\Omega, \mathcal{F}, \mathbb{P})$.} \vEE{In the following, for a given stochastically continuous rf $U(t),t\in \TTT$ we shall suppose that it is also separable and jointly measurable. In view of \cite{Doob} a separable and jointly measurable version of a stochastically continuous rf  exists.
}

\begin{theorem}\label{th1}
	If  $Z(t), t\in \R^d$ is a  non-negative  \bE{stochastically continuous}  rf such that \rrd{\eqref{MS}},
	 \eqref{norm} and \eqref{shi} hold,  then   the  constants defined in \eqref{pcf} are finite  for all $\delta \ge 0$ and further   $\lim_{\delta\downarrow 0} H_{Z}^{\delta }= H_{Z}^0$.
\end{theorem}
In the following
$$ S_\delta (Z)= \int_{t\in \delta \Z^d }  {Z(t)}  \lambda_\delta(dt), \quad \delta\ge 0,
$$
where  $\lambda_0(dt)=\lambda(dt)$ is the Lebesgue measure on $\R^d$  and for $\delta>0$ $\lambda_\delta(dt)= \delta ^d \lambda(dt)$ with $\lambda(dt)$  the counting measure on $\delta \Z^d$. \\
\cEE{
Clearly, $S_\delta(Z)$ is a random variable (rv) for any $\delta>0$.
If $\delta=0$  (recall $0\Z^d$ denotes simply $\R^d$), since we consider  $Z$ to be jointly measurable,
then \kkk{supposing that \eqref{norm} holds and using} that $Z$ is non-negative
(recall we assume that the probability space is complete)
it follows that $S_0(Z)$ is a rv; see \cite{Doob}[Thm 2.7, 2.8] for details in case $d=1$.}

\begin{corollary}
Suppose that  $Z$ satisfies the assumptions  of \netheo{th1} and further $\pk{\sup_{t\inr^d} Z(t)>0}=1$. If
$ S_\eta(Z)  = \IF$ almost surely for some $\eta\ge 0$, then $H_{Z}^\delta=0$ for all $\delta \ge 0$. Conversely,
if $H_{Z}^\delta=0$ for some $\delta \ge 0$, then $S_\eta(Z)=\IF$ almost surely for all $\eta\ge 0$.
\label{kor1}
\end{corollary}
	\BRM\label{rem1}
\begin{enumerate} [i)]
	\item \label{rem1:A}
 The assumption    $\pk{\sup_{t\inr^d} Z(t)>0}=1$ in \nekorr{kor1} cannot be removed.
 Taking for instance   $Z(t)= V/p,t\inr^d$ with $V$ a Bernoulli rv  such that $\pk{V=1} =p\in (0,1)$,
 we get  $H_{Z}^\delta=0$ for all $\delta\ge 0$. However $S_1(Z)=0$ with probability \rrd{$1-p>0$}.
	\item \label{rem1:B} Let $Z(t),t\inr^d$ be a stationary rf satisfying the assumptions of \netheo{th1}.
By the definition and the fact that $Z$ and $\tilde Z= 1+ Z$  are both stationary, the corresponding Pickands constants exist and we simply have  $H_{Z}^\delta=H_{\tilde Z}^\delta$ for all $\delta\ge 0$.  Clearly $S_{1}(\tilde Z)=\IF$ and $\pk{\sup_{t\inr^d} \tilde Z(t)>0}=1$, hence \nekorr{kor1} implies $H_Z^\delta=0$ and moreover $S_\eta(Z)=\IF, \eta \ge 0$ if further  $\pk{\sup_{t\inr^d}  Z(t)>0}=1$.
%Note that as shown in \ref{rem1:A}  above, in general  the stationarity of  $Z $ does not imply $S_1(Z)=\IF$ almost surely.
\end{enumerate}
	\ERM
In the rest of this section we consider $Y(t),t\inr^d$ a stationary max-stable rf
as in the Introduction with spectral rf $Z$ and de Haan representation \eqref{eq1}.
Suppose next that $Z$ has sample paths almost surely in the space $D=D(\R^d, [0,\IF))$ of generalised  \cadlag\ functions  $f: \R^d \to [0,\IF)$ equipped with the  $\sigma$-field $\mathcal{D}=\sigma(\pi_t, t\in \TT_0)$ generated by the projection maps $\pi_t: \pi_t f= f(t), f\in D$ and
$\TT_0$ a dense subset of $\R^d$.  See e.g., \cite{Svante,MartinE} for the definition and properties of generalised  c\`{a}dl\`{a}g functions.
In view of \cite{Htilt}[Thm 6.9] (take $\alpha=1$ and $L=B^{-1}$ therein)  the stationarity of  $Y$ is equivalent with the validity of
\bqn{
	\label{rinashero} \Ex{ Z(h) F(Z)} = \Ex{ Z(0) F(B^h Z)}
}
 for all $h\in \R^d$ and all  0-homogeneous measurable functionals $F: D \mapsto [0,\IF)$ ($0$-homogeneous
 means $F(cf)=F(f), \forall c>0, f\in D$) with $B^h Z(\cdot)= Z(\cdot-h), h\in \R^d$.
 For discrete max-stable processes, \eqref{rinashero}  is stated in \cite{Hrovje}[Eq.\ (5.2)],  see also \cite{PH2020} for other equivalent formulations. We note in passing that \eqref{rinashero} is initially derived  for  $Z$ as in Example \ref{ex1} below  in \cite{DM}[Lem 5.2].

Clearly, if $Z$ is stationary, then $Y$ is stationary too. However this instance is not interesting,
since as shown in Remark \ref{rem1}, \ref{rem1:B}  $H_{Z}^0=0$ \kkk{in this case}.
We discuss  below  two other examples such that max-stable rf $Y$ is stationary.

\begin{example}\label{ex1}
It is known from \cite{K2010,Htilt} that if $X(t)= W(t)- Var(W(t))/2, t\in \R^d$ with $W$
a centered Gaussian  rf with stationary increments and almost surely continuous sample paths,
then both \eqref{MS} and \eqref{shi} hold with $Z(t)=e^{X(t)}$. Using for instance \cite{Vitale1}[Thm 1] we have that \eqref{norm} holds. Moreover, in view of \cite{kab2009} the corresponding max-stable  rf $Y$ is stationary.
\end{example}

\begin{example}\label{ex3}
Let $L(t),t\in \R^d$ be a non-negative deterministic measurable function such that $\int_{\R^d} {L (t)}\lambda(dt) = 1, $
with $\lambda(dt)$ the Lebesgue measure on $\R^d$. It follows easily that
$ Z(t) =   L(t-N) / p(N) ,t\inr^d$ with $N$
an $\R^d$-valued rv having a  positive density function $p(t)>0, t\in \R^d$ satisfies \eqref{shi}.
Note that for this case condition \eqref{norm} reads
\bqn{ \label{reads}
\int_{x\in \R^d}  \sup_{t\in \rrd{[0,T]^d}}  {L(t-x)}  \lambda(dx)< \IF, \quad \forall\, T>0.
}
\end{example}

As shown in \cite{HBernulli}[Thm 1, Eq.\ (3.5)]  without any further assumption on the max-stable stationary rf $Y$,  \tE{for $\delta=\eta> 0$} we have
	\bqn{ 	\label{drejt}
		H_{Z}^\delta = \Ex*{ {Z(0)}   \frac{\sup_{t\in \delta \Z^d }  {Z(t)}  }{S_\eta(Z) }   } \in [0,\IF).
	}
In the next result we show that the above holds under some weak assumptions also for $\delta=0$. It turns out that under \eqref{laps} below  we can obtain also  the first formula in \eqref{DMY}.

\begin{proposition} Let   $Z(t),t\inr^d$ be a non-negative  rf with almost surely sample paths in $D$. If  \eqref{MS}, \eqref{norm}
	and  \eqref{rinashero} hold and  $H_{Z}^0>0$, then \eqref{drejt} holds for $\delta=\eta= 0$.
Moreover, \eqref{drejt} holds also for $\delta=0, \eta>0$ or $\delta>0, \eta=k\delta, k \inn$ if further
\bqn{ \label{laps}
 \{ \vEE{S_0(Z)} < \IF \} \subset   \{ S_\eta(B^r Z) \in (0,\IF)   \}, \ \ \forall r\in \delta \Z^d
}
almost surely.
\label{TM2}
\end{proposition}

\begin{remark} \label{REMA}  	
%\begin{enumerate} [i)]
%	\item  \label{REMA:A}
	\vEE{It is shown in the proof of \neprop{TM2} that \eqref{laps}  is implied by the assumption
$\pk{Z(0)>0}=1$, which is satisfied for the choice of $Z$ as in Example 1. In particular, \neprop{TM2} extends
\cite{SBK}[Thm 2 and Thm 3] and \cite{DiekerY}[Prop 2].}
Moreover, as  discussed in Example \ref{ex.3} below, condition \eqref{laps} cannot be removed.
\end{remark}
\begin{example}\label{ex.3}
Let $L(t),t\in \R^d$ be as in Example \ref{ex3} and suppose further that $L(t)>0, t\inr^d$,   \eqref{reads} holds and $L \in D(\R^d, [0,\IF))$.  Then the claim of \netheo{th1} follows. Moreover, \neprop{TM2} implies for all $\eta >0$  %This can be shown also directly since  \eqref{vivet} implies
\bqn{ \label{ishk}
	H_{Z}^0=
\int_{\R^d}   L (s)
	\frac{ \sup_{t\in \R^d}  {L(t+s)}  } { \eta^d \sum_{t\in \eta \Z^d}  {L(t+s)} } \lambda(ds)=
	 \sup_{t\in \R^d} {L (t) } ,
}
\tE{where the last equality follows by  (\ref{drejt}) with $\delta=0$ (recall $\int_{\R^d} L(t) \lambda(dt)=1$).}
It is  mentioned in \cite{DiekerY} that for  $d=1$  and $L$  the standard Gaussian density on $\R$,
the above identity could be verified numerically. As shown by Dmitry Zaporozhets (personal communications)
the last two equalities in \eqref{ishk} can be derived  \vEE{utilising  the translation invariance of both the
Lebesgue measure and the counting measure, respectively} and applying the Fubini-Tonelli theorem. Specifically, setting  for simplicity below $\eta=1$
\bqny{
	\int_{\R^d}
	\frac{  L (s) } { \sum_{t\in \Z^d}  {L(t+s)} } \lambda(ds) &=& %\sup_{t\in \R^d} {L (t) }
		\sum_{ r \in \Z^d}   \int_{ r+ [0,1)^d }
	\frac{ L(s)  } {   \sum_{t\in   \Z^d}  {L(t+s)} } \lambda(ds)\\
	 &=& 	\sum_{ r \in \Z^d}   \int_{ [0,1)^d }
	\frac{ L(r+s)  } { \sum_{t\in   \Z^d}  {L(r+t+s)} } \lambda(ds)\\
	 &=& 	   \int_{ [0,1)^d }\sum_{ r \in \Z^d}
\frac{ L(r+s)  } { \sum_{t\in   \Z^d}  {L(t+s)} } \lambda(ds)=1.
}
Our initial proof of  \neprop{TM2} was asymptotic in nature. A modification   of Dmitry's arguments and \eqref{rinashero} have led to the current proof of \neprop{TM2}.\\
Note that if we take for instance $L(t)=\ind{t\in [0,1]},  d=1$ with $\ind{\cdot}$ the indicator function, then  $H_Z^0$ cannot be given by \eqref{drejt} with  $\eta>2$, since  the first equality in \eqref{ishk} gives $H_Z^0=\IF$, which is a contradiction to the fact that $H_Z^ 0< \IF$. For this choice of $\eta$ we have that condition \eqref{laps} is not satisfied.
\end{example}

\def\HFA{H_{\abs{f},m}^\delta }

\section{Two Extensions}
Motivated by \cite{Genna04c, Genna04} we consider first  an alternative definition of Pickands constants $\HFA$ defined with respect to a collection of functions $f_s: E \to \R $ with $(E, \mathcal{E})$ a measurable space equipped with some $\sigma$-finite  measure $m$. In the second part of this section we discuss the constant $\mathcal{H}_{Z}^\delta$ defined in \eqref{hexh} proving also its continuity at $\delta=0$.

\subsection{Behaviour  of $\HFA$ at $\delta=0$} Let $m$ be a  $\sigma$-finite measure  on some measurable space $(E, \mathcal{E})$ and
$f={f_s(z),s\in \R^d, z\in E}$ real functions with $f_s \in \mathcal{L}^1 (m)$  ($\mathcal{L}^\alpha(m)$ is the set of all functions $g$ such that $\int_{E}|g(x)|^\alpha m(dx)<\infty$ and $\alpha>0$).   We define next for any $\delta =1/n,n \inn$ or $\delta=0$
$$ \HFA = \limit{T} \frac{1}{T^d}
\int_{E }  \sup_{t\in [0,T]^d \cap \delta \Z^d} \abs{f_t (z)} m(dz),
$$
where  $0\Z^d$ equals $\mathbb{Q}^d_*$ with $\mathbb{Q}_*$ the set of dyadic rational numbers
\rrd{\kkk{\{}$\frac{k}{2^n}$: $k\in \Z$, $n\inn$\kkk{\}}}.
 Clearly, the above limit is in general not defined. In order to include also the case $\delta=0$ we shall assume further that
\bqn{ \int_E \sup_{t\in \KK \cap \mathbb{Q}_*^d}  \abs{f_t(z)}   m(dz) < \IF
\label{norm2}
  }
and  similarly to \eqref{shi}
\bqn{ \int_E \sup_{t\in K \cap \mathbb{Q}_*^d }  \abs{f_{t+c}(z)}   m(dz) = \int_E \sup_{t\in K \cap \mathbb{Q}_*^d }  \abs{f_{t}(z)}   m(dz), \quad
	\forall c\in \mathbb{Q}_*^d
\label{genA}
}
hold  for all  compact sets $\KK \subset \R^d$.

If $L$ is as in Example 2, then taking $f_t(z)= L(z-t), t, z\inr^d
$,  $E=\R^d$ equipped with Borel $\sigma$-field and $m(dz)$ the Lebesgue measure on $\R^d$ we have that
 \eqref{norm2} holds and further \eqref{genA} is satisfied since $m(dz)$ is shift-invariant.

\begin{example}\label{ex5}
Let $Y(t),t\inr^d$
be  a symmetric $\alpha$-stable  stationary  rf with locally bounded sample paths, $\alpha \in (0,2)$ and representation
$$ Y(t) =  \int_E f_t(z) M(dz) , \quad t\in \R^d ,
   $$
where $M$ is a symmetric $\alpha$-stable random measure on $E$ with control measure $m$ and $f_t\in \mathcal{L}^\alpha(m), t\inr^d$, see \cite{Genna04c}. It follows that both \eqref{norm2} and \eqref{genA} hold, see \cite{Genna04c} for the case $d=1$ and \cite{Roy1} for the case $d>1$.
\end{example}
With the same arguments as in the proof of \netheo{th1}, we have that
\kkk{under} \eqref{norm2} and \eqref{genA}, \kkk{constant} $\HFA$ is finite,  non-negative and further
\bqn{ \label{dichte}
	\limit{n} H_{\abs{f}, m}^{2^{-n}} =
H_{\abs{f}, m}^{0} .
}

\BRM For $Y$ as in Example 4 explicit formulas for $H_{\abs{f}, m}^{\delta} $ are derived in \cite{Genna04c,Genna04,MR2384479}. Utilising the relation between $\alpha$-stable and max-stable processes, see \cite{kab2009a}, the aforementioned formulas imply \eqref{drejt} when $\delta=\eta\ge 0$.
\label{remGenna}
\ERM

\subsection{Continuity of $\mathcal{H}_Z^\delta$ at $\delta=0$}
Let $W_z(t),t\inr^d, z\in [0,1] $ be a  centred  Gaussian  rf with stationary increments, almost surely continuous sample paths  and variance function $\sigma^2_z,z\in [0,1]$ such that $\sigma^2_z(0)=0$ for all $z\in [0,1].$  We formulate below the assumptions on the variance functions $\sigma_z^2$ imposed in \netheo{th11} below. Specifically, we shall assume that
\bqn{ \label{sz}
\lim_{w\to z}\sigma_w(t)=\sigma_z(t), \quad \forall z \in [0,1], \forall t\inr^d
}
and further for  some $C_{0},C_{\infty}$ positive and $\nu_{0},\nu_{\infty}\in (0,2] $
	\bqn{ \label{kN}
		\limsup_{\norm{t}\to 0}	 \frac{ \sigma_z^2(t)}{\norm{t}^{\nu_{0}}} \le C_{0}\, \ {\rm and}\, \
		\limsup_{\norm{t}\to \infty}	 \frac{ \sigma_z^2(t)}{\norm{t}^{\nu_{\infty}}} \le C_{\infty}
	}
	hold  for all $z\in [0,1]$, \rrd{where $\norm{\cdot}$ denotes the Euclidean norm on $\R^d$.}

\begin{theorem}
	If \eqref{sz} and  \eqref{kN} are satisfied, then
	for any $\delta\geq 0$ and $Z_z(t)= \exp( W_z(t)- \sigma_z^2(t)/2)$
	\bqn{
		\mathcal{H}_{{{Z}}}^\delta =   \limit{T} T^{-d} \int_0^1 \Ex*{\sup_{ t\in [0,T]^d \cap \delta \Z^d }Z_z (t)  } dz=\int_0^1  H_{Z_z}^\delta dz \in [0,\infty)
	}
	and furthermore $\lim_{\delta \downarrow 0} \mathcal{H}_{Z}^\delta= \mathcal{H}_{Z}^0$.
	\label{th11}
\end{theorem}

\BRM \begin{enumerate} [i)]
	\item
 A sufficient condition for $H_{Z_z}^\delta, z\in [0,1]$ to be positive is
\bqn{\label{KD}
	\limit{\norm{t}} \frac{ \sigma_z^2(t)}{ \ln \norm{t}} > 8 d,
}
see \cite{debicki2017approximation}.
If the above holds for $z\in [0,1]$ on a set with non-zero Lebesgue measure, then under the assumptions of \netheo{th11} we have that
$\mathcal{H}_{Z}^\delta>0$ for any $\delta \ge 0$.
\item  In order to define  rf's $W_z,z\in [0,1]$ with stationary increments, we need to determine variance function $\sigma^2_z(t)= Var( W_z(t+s)- W_z(t)), s,t\in \R^d$, which is a negatively (or conditionally negatively) definite function.  One  instance is to start with a negatively definite variance function $\sigma^2(t), t\inr^d$ and then define $\sigma_z(t)= q(z) \sigma(t)$ for some  continuous function $q$ not identical to zero. Such $\sigma_z$ clearly satisfies \eqref{sz} and if $\sigma^2(t)\le C \norm{t}^\nu, t\inr^d$,  then \eqref{kN} also holds.
Another instance is to \kkk{consider} $\sigma_z$ as in \cite{MR2206343}, namely
\bqn{\label{CL}
	\sigma^2_z(t) = \norm{t}^\lambda r_z(t/\norm{t}), \quad \lambda \in (0,2], t\inr^d,
}
where $r_z$ is a non-negative function defined on  the unit sphere on $\R^d$
determined by the norm $\norm{\cdot}$.
\end{enumerate}
\ERM

\def\ve{\varepsilon}
\section{Proofs}
\prooftheo{th1}
  In view of \eqref{norm} and
 \eqref{shi}
\begin{equation}\label{sonst}
	\Ex*{ \sup_{t\in \KK  \cap \delta   \Z^d} Z(t+c)}= \Ex*{ \sup_{t\in \KK \cap \delta   \Z^d} Z(t)}< \IF
\end{equation}holds for all
$c\inr^d, \delta  \ge 0$ and all
  compact sets $\KK\subset   \R^d$.
Next, assume for notational simplicity that
$d=1$.
In the light of \eqref{sonst}  $\Ex{ \sup_{t\in [0,T]  \cap \delta \Z } {Z(t)}}, T>0$ is
a subadditive function of $T$ for any $\delta\ge 0$ implying that  $H_Z^\delta$ is a finite non-negative constant (recall $Z$ is non-negative)
in view of Fekete's lemma.
Further, to show the second part of the theorem it is enough to take $T= a_\delta n$ where $a_\delta=\floor{1/\delta} \delta $ with
$\floor{x}$  the integer part of $x>0$ and $n \inn,\,\delta \in (0,1/2)$. Utilising \eqref{sonst} to derive the last equality in the next calculations
we have
  \begin{eqnarray}\lefteqn{  \Ex*{ \sup_{t\in [0,T]}   {Z(t)}   }-\Ex*{ \max_{t\in [0,T] \cap \delta \Z}  {Z(t)} }}\nonumber\\
  	&=&  \Ex*{ \max_{ 1 \le i \le n} \sup_{t\in [a_\delta(i-1),a_\delta i]}   {Z(t)}  -  \max_{ 1 \le i \le n}  \max_{t\in [a_\delta(i-1),a_\delta i] \cap \delta  \Z}  {Z(t)} }\nonumber\\
&\leq &  \Ex*{ \max_{ 1 \le i \le n} \Bigl(\sup_{t\in [a_\delta(i-1),a_\delta i]}   {Z(t)}  -  \max_{t\in [a_\delta(i-1),a_\delta i] \cap \delta  \Z}  {Z(t)} \Bigr) }\nonumber\\
	&\le &
	  \sum_{i=1}^n \Ex*{ \sup_{t\in [a_\delta(i-1),a_\delta i]}  {Z(t)}  -  \max_{t\in [a_\delta(i-1),a_\delta i] \cap \delta  \Z}  {Z(t)} }\nonumber\\
	&= & n \Ex*{ \sup_{t\in [0,a_\delta]} {Z(t)} -  \max_{t\in [0,a_\delta] \cap \delta  \Z}  {Z(t)} }.\label{ineqcont}
\end{eqnarray}
The assumption that $Z$ is stochastically continuous implies that any
 dense subset of $\R^d$ is a separant for the separable  rf $Z$ (see e.g., \cite{MR2541742}[Thm 2.9]) and therefore since $\lim_{\delta \downarrow 0} a_\delta=1$ the following convergence in probability
% https://encyclopediaofmath.org/wiki/Separable_process
\bqny{
	\max_{t\in [0,a_\delta] \cap  \delta \Z}   {Z(t)}  \to   \sup_{t\in [0,1]}  {Z(t)}
}
 holds as $\delta \downarrow 0$. Hence by \eqref{norm}
 \rrd{and the dominated convergence theorem} %? Below we use mononicity but above we cannot use?}
\bqn{\label{tog} \lim_{\delta \downarrow 0 } \Ex*{ \max_{t\in [0,a_\delta ] \cap \delta  \Z}  {Z(t)} }
	% =\lim_{\delta \downarrow 0 } \Ex*{ \sup_{t\in [0,a_\delta ] }  {Z(t)} } =
	=\Ex*{ \sup_{t\in [0,1]} {Z(t)}}< \IF.
}
Consequently   \eqref{ineqcont} implies
\bqny{
	 \frac{1}{T} \Ex*{ \max_{t\in [0,T] \cap \delta  \Z}  {Z(t)}} & \leq &
 \frac{1}{T}	 \Ex*{ \sup_{t\in [0,T]}   {Z(t)}}\\
	 &\leq  & a_\delta^{-1} \Ex*{ \sup_{t\in [0,a_\delta]} {Z(t)} -  \max_{t\in [0,a_\delta] \cap \delta \Z}  {Z(t)} }	
	+   \frac{1}{T}\Ex*{ \max_{t\in [0,T] \cap \delta  \Z}  {Z(t)}}.
}
Since $H^\delta_Z$ exists and is finite for all $\delta \ge 0$ we have
$$  \limit{T}  \frac{1}{T}	 \Ex*{ \sup_{t\in [0,T]}   {Z(t)}} =
 \limit{n}  \frac{1}{ n a_\delta }	 \Ex*{ \sup_{t\in [0,n a_\delta]}   {Z(t)}} = H^0_Z$$
 and
$$  \limit{T}  \frac{1}{T}	 \Ex*{ \sup_{t\in [0,T] \cap \delta \Z }   {Z(t)}} =
\limit{n}  \frac{1}{ n a_\delta }	 \Ex*{ \sup_{t\in [0,n a_\delta] \cap \delta Z}   {Z(t)}} = H^\delta _Z  \le  H^0_Z<\IF.$$
Thus letting \rrd{$n$} to infinity we obtain from the above inequalities
% for all $\delta>0$ sufficiently small and all $C>1$
%since the Pickands constants exist
\bqn{
0 & \leq &
	H_{Z}^{0 }-H_{Z}^{\delta } \leq   \rrd{a_\delta^{-1}}\Ex*{ \sup_{t\in [0,a_\delta]} {Z(t)} -  \max_{t\in [0,a_\delta] \cap \delta \Z}  {Z(t)} },
%	\le
%	C\Ex*{ \sup_{t\in [0,1]} {Z(t)} -  \max_{t\in [0,1] \cap \delta \Z}  {Z(t)}},		
}
which together with \eqref{tog} yields   $\lim_{\delta \downarrow 0} H_{Z}^{\delta }= H_{Z}^0$.
\QED

\proofkorr{kor1}
First, let us notice that $H_{Z^{}}^\delta =0$ for some $\delta\geq 0$ if and only if $S_\delta(Z)=\IF$ almost surely, which is a direct implication of \cite{dom2016,kab2009a} and \cite{Genna04c,Genna04,MR2384479}; this has been already discussed in \cite{MR2453345,debicki2017approximation, PH2020} for the case $d=1$ and \cite{HBernulli} for the $d$-dimensional discrete setup. Note in passing that in \cite{dom2016} $Z$ is such that $\pk{\sup_{t\in \R^d} Z(t)> 0}=1$.   So if for some $\eta>0$ we have  $H_{Z}^\eta=0$,   then
$S_{\eta/k}( Z)  = \IF$ almost surely for any $k\in N$ and  $H_{Z }^{\eta/k}=0$. Consequently,  \netheo{th1} implies that $\limit{k}H_{Z }^{\eta/k}=H_{Z }^0=0$ and hence $0\le H_{Z}^\eta
\le H_{Z }^{0}= 0$ for all $\eta\ge 0$. Conversely, if for some $\eta\ge 0$ we have  $S_{\eta}( Z)  = \IF$, then from the above equivalence
$H_{Z}^\eta=0$ for all $\eta\ge 0$, hence the proof is complete.
\QED

\proofprop{TM2}
The claim for $\delta=\eta=0$ follows as in \cite{debicki2017approximation} where $d=1$ is considered and its proof is therefore omitted.\\
% \cEE{An alternative proof is to use the approach developed for $\alpha$-stable processes, see Remark \ref{remGenna}.}\\
 Let next $D=D(\R^d, [0,\IF))$ be the space of generalised \cadlag\  functions $f: \R^d \mapsto [0,\IF)$, which can be  equipped with a metric that turns it into a Polish space. The corresponding Borel $\sigma$-field in $D$  denoted by $\mathcal{D}$ agrees with the $\sigma$-field $\sigma(\pi_t, t\in \TT_0)$ for any $\TT_0$ a dense subset of $\R^d$, see e.g., \cite{MartinE}[Thm 7.1].   Denote by $\Theta(t),t\inr^d$ a  rf with almost surely sample paths in $D$ defined by
$$\pk{ \Theta \in A}= \Ex*{ Z(0) \mathbb{I}( Z/Z(0)\in A)}/\Ex*{ Z(0)}, \quad \forall A \in \mathcal{D},$$
where we interpret $0:0$ as 0 and set   $\mathbb{I}( x\in A)$ equals 1 or 0 if $x\in A$ or $x\not \in A$, respectively.
\vEE{Note that $\Ex*{ Z(0)}$ equals 1, but we leave it since the same formula is applied to $Z_*$ below}. Note further that since $D$ is Polish, by  \cite{Vara}[Lemma on p. 1276] we can realise both $Z$ and $\Theta$ in the same complete non-atomic probability space which we assume for notational simplicity below. By \eqref{rinashero} for any measurable functional $F: D \mapsto [0,\IF]$ which is $0$-homogeneous
\bqn{
	\label{rinashero2} \Ex{ \Theta (h) F(\Theta)} =
	\Ex{  \mathbb{I}( Z(0)\not=0)  Z (h)  F(Z)}=
	 \Ex{  \mathbb{I}(\Theta(-h) \not=0 ) F(B^h \Theta) } .
}
In view of \nekorr{kor1} the assumption $H^0_{Z}>0$ implies $\pk{S_\eta(Z)< \IF}=p_\eta >0$  for all $\eta\ge 0$.
 If $p_0=1$, then  also \cL{ $\pk{S_0(\Theta) < \IF}=1$} follows and thus
 \bqn{ \{ S_\eta(\Theta) < \IF \} \subset \{ S_0(\Theta)< \IF\}
 \label{shtu}
}
almost surely. \\
   Assume next that $p_0\in (0,1)$. The rf   \kdd{$Z_*(t)=Z(t) \vEE{\lvert  S_0(Z) \rrd{=} \IF}, \,t
 	\in \R^d $} has almost surely sample paths in $D$ and satisfies \eqref{rinashero}.  Denote by $\Theta_*$ the corresponding  rf of $Z_*$ defined by the change of measure as above. We have that $\Theta_*$ has the same law as $\Theta \lvert S_0(\Theta)= \IF$. 	Since $
 S_0(  Z_*) =\IF$ with probability 1,
 applying \nekorr{kor1}, we obtain
$ S_\eta(Z_*)=\IF$ almost surely and thus $S_\eta  (\Theta_*)=\IF$ almost surely  i.e.,
$$\pk{ S_\eta(\Theta)=\IF, S_0(\Theta)= \IF}/\pk{S_0(\Theta)= \IF}=1$$
 \rrd{(implying $ \{ S_0(\Theta) = \IF \} \subset \{ S_\eta(\Theta) = \IF\}$ almost surely)}, hence  \eqref{shtu}  holds. If $\pk{S_0(Z)< \IF}=0$ then by \nekorr{kor1}
 $\pk{S_\eta(Z)< \IF}=0$ for all $\eta>0$, which in turn implies $\pk{ S_\eta(\Theta) < \IF} =0$. The assumption $p_0=0$ implies also
 $\pk{ S_0(\Theta) < \IF} =0$, hence again \eqref{shtu} holds.\\
 Next, suppose that  $\pk{S_0(\Theta)< \IF}>0$. Since then $\pk{S_0(Z)< \IF}>0$, we can define as above \kdd{$Z_*(t)=Z(t) \vEE{\lvert S_0(Z) <  \IF},\, t\in \R^d $}.
\COM{\bqny{ \pk{\Theta_*\in A} &=& \Ex{ Z_*(0)\ind{  Z_*/Z_*(0) \in A } } /\Ex{Z_*(0)} \\
	&=& \Ex{Z(0) \ind{ S(Z)< \IF } \ind{  Z \ind{ S(Z)< \IF } \in A}   / ( Z(0)\ind{ S(Z)< \IF }) \in A } /\Ex{ Z \ind{ S(Z)< \IF }}\\
	&=&  \Ex{\ind{ S(\Theta) < \IF, \Theta \in A}}/\Ex{ Z \ind{ S(Z)< \IF }}= \pk{ \Theta \in A \lvert S(\Theta) < \IF}
}
}
The corresponding $\Theta_*$ has the same law as $\Theta \lvert S_0(\Theta)< \IF$. Since $\pk{S_0(Z_*)< \IF}=\pk{S_0(\Theta_*)< \IF}=1$, by \cite{PH2020}[Thm 2.8]
 $$\pk*{\lim_{\norm{t}\to \IF, t\in \eta \Z^d} \Theta_*(t) =0}=1,$$
  with $\norm{\cdot}$ some norm on $\R^d$. The latter is equivalent with
 $\pk{S_\eta(\Theta_*)< \IF}=1$ see \cite{HBernulli}[Condition A2, A4].
As in the proof of \cite{HBernulli}[Lem A.2]
 $\pk{ S_0(\Theta) < \IF} =0 \iff p_0=0$ and thus   the reverse inclusion to \eqref{shtu} holds implying
\begin{equation}\label{LUK}
	\{  S_\eta(\Theta) < \IF \} =\{  S_0(\Theta)< \IF\}, \ \forall \eta>0
\end{equation}	
   almost surely. Since further  $\pk{S_\eta(\Theta)>0}=1$ for all $\eta\ge 0$, which follows from
$\pk{\Theta(0)=1}=1$ and    the fact that  $\Theta$ has \vEE{paths in $D$ almost surely}, \vEE{we have}
almost surely for all $\delta,\eta \in [0,\IF)$
% for $\delta=0$ and all $\eta>0$  and for  $\delta>0$ and all $\eta= k\delta$, $k\inn$    (recall that $\pk{ \Theta(0)=1}=1$)
\bqn{ \label{lime2}
	 \frac{ 1 }{S_\eta(\Theta)} =\frac{ 1 }{S_\eta(\Theta)} \frac{ S_\delta(\Theta)} {S_\delta(\Theta)}\Theta(0).
}	
  Set next   $S_\eta(f)= \int_{ \eta \Z^d} f  (t)\lambda_\eta(dt), \eta\ge 0$ and
  $$M_\eta(f)= \sup_{t\in \eta \Z^d} \tE{f }(t),  \ \    M_0(f)=  \sup_{t\in \TT_0} f(t), f\in \vEE{D}, \eta >0.
  $$
  %	 with $\TT_0$ a dense countable set of $\R^d$ which is a separant for the separable process $Z$.
   Both maps
   $S_\eta(\cdot)$ and $M_\eta(\cdot), \eta \ge 0$ are measurable and by the separability of $Z$ we have that   $M_0(Z)$ has the same law as $\sup_{t\in \R^d} Z(t)$. Hence using further the definition of $\Theta$ and \eqref{lime2}, by the Fubini-Tonelli theorem  for all $\delta\not=\eta, \delta \ge 0, \eta>0$ we have
  \bqn{ \lefteqn{\Ex*{ Z (0)  \frac{ \sup_{t\in  \delta \Z^d} Z (t)}{S_{\eta }(Z)}}}\notag\\
	&=& \Ex*{  \frac{ \sup_{t\in  \delta \Z^d} \Theta (t)}{S_{\eta }(\Theta)}} \notag\\
	&=& \Ex*{  \frac{ \sup_{t\in  \delta \Z^d} \Theta (t)}{S_{\eta }(\Theta)} \frac{ S_\delta(\Theta)}{ S_\delta(\Theta)} \Theta(0)} \notag\\
	&=&  \int_{\delta \Z^d}\Ex*{  \frac{ M_\delta (  \Theta)}{ S_\delta(\Theta)}\frac{ \Theta (0)}{ S_{\eta }(\Theta) }  \Theta  (s) } \lambda_\delta (ds) \notag\\
	&=&  \sum_{i \in \eta \Z^d}	 \int_{r\in [0,\eta)^d \cap \delta \Z^d }\Ex*{  \frac{ M_ \delta (\Theta) }{ S_\delta(\Theta) }\frac{ \Theta (0)}{ S_{\eta }(\Theta)}  \Theta (i+r) } \lambda_\delta (dr),
	 \label{i2}
}
\vEE{where in the last equality we used also the translation invariance of $\lambda_\delta$.}
{Now by \eqref{rinashero2} for all  $i,r \in \R^d$ %, r\in \delta \Z^d$
% we have
\begin{equation}\label{schlb}
	\Ex*{   \Theta (i+r) \frac{M_ \delta (\Theta) }{ S_\delta(\Theta) }\frac{ \Theta (0)}{ S_{\eta }(\Theta)}}
   = \Ex*{ \frac{\ M_ \delta (B^{r+i}\Theta)}{ S_\delta(B^{r+i}\Theta)}     \frac{\Theta (-i-r) }{S_{\eta }(B^{r+i}\Theta)}   }.
\end{equation}
\vEE{Note in passing that  $\Theta(-i-r) \mathbb{I}(\Theta(-i-r) \not=0 )=\Theta(-i-r) $ almost surely.} \\
 If  $\delta=0, \eta>0$  or  $\delta>0$ and $\eta= k\delta$, $k\inn$,
   then almost surely for all $i \in \eta \Z^d, r\in \delta \Z^d$
\begin{equation}\label{ii}
	 M_ \delta (B^{r+i}\Theta)= M_ \delta (\Theta),  \  S_\delta(B^{r+i}\Theta)= S_\delta(\Theta), \  S_{\eta }(B^{r+i}\Theta)
   =S_{\eta }(B^{r}\Theta)
\end{equation}
and from \eqref{laps}%, \eqref{LUK} %we have further
%where  the last claim in \eqref{ii} follows from the assumption that
\bqny{   \pk{S_\delta (\Theta) < \IF} %&=&\pk{S (\Theta) < \IF}\\
     &=& \vEE{\pk{S_0 (\Theta) < \IF}}\\
	&=& \vEE{\Ex{ Z(0) \ind{  S_0(Z)< \IF}}}\\
	&=& \Ex{ Z(0) \ind{  \vEE{S_0(Z)}< \IF,  S_\eta(B^r Z) \in (0,\IF) }}\\
	&=& \Ex{ Z(0) \ind{  \vEE{S_{\delta}}(Z)< \IF,  S_\eta(B^r Z) \in (0,\IF) }}\\
	&=& \Ex{ Z(0) \ind{  \vEE{S_{\delta}}(Z/Z(0))< \IF,  S_\eta(B^r Z/Z(0)) \in (0,\IF) }}\\
	&=& \pk{
		 S_\delta(\Theta) < \IF, S_\eta(B^r \Theta) \in (0,\IF)}, \ \forall r\in \delta \Z^d,
}
\COM{

\bqny{ 1&=&\Ex{Z(0)}=\Ex{ Z(0) \ind{ S_\eta(B^r Z) \in (0,\IF), S(Z)< \IF}}\\
	&=&  \Ex{ Z(0) \ind{ S_\eta(B^r Z/Z(0)) \in (0,\IF), S(Z/ Z(0)) < \IF}} \\
	&=& \pk{
	S_\eta(B^r \Theta) \in (0,\IF), S(\Theta) < \IF}\\
	&=& \pk{
	S_\eta(B^r \Theta) \in (0,\IF), S_\delta(\Theta) < \IF}, \ \forall r\in \delta \Z^d.
}
}
\vEE{where the third last line above follows from \eqref{LUK}}. Consequently,  almost surely
\bqn{\label{LLI}  \frac {S_{\eta }(B^r\Theta) }{S_{\eta }(B^r\Theta)}  \ind{S_\delta(\Theta) <\IF}    =  \ind{S_\delta(\Theta) <\IF} .
}
Hence, for these choices of $\delta $ and $\eta$,   \eqref{i2}-\eqref{LLI} and
 the Fubini-Tonelli theorem yield
  \bqny{ \Ex*{ Z(0)  \frac{ \sup_{t\in  \delta \Z^d} Z(t)}{S_{\eta }(Z)}}
&=&  	 \int_{r\in [0,\eta)^d \cap \delta \Z^d }\Ex*{ \frac{\ M_ \delta (\Theta)}{ S_\delta(\Theta)}  \sum_{i \in \eta \Z^d}  \frac{\Theta (-r-i) }{S_{\eta }(B^r\Theta)}   } \lambda_\delta (dr) \\
&=&  	 \eta^{-d} \int_{r\in [0,\eta)^d \cap \delta \Z^d }\Ex*{ \frac{\ M_ \delta (\Theta)}{ S_\delta(\Theta)}   \frac {S_{\eta }(B^r\Theta) }{S_{\eta }(B^r\Theta)}  \ind{S_\delta(\Theta) <\IF} } \lambda_\delta (dr) \\
&=&  \eta^{-d} 	 \int_{r\in [0,\eta)^d \cap \delta \Z^d }\Ex*{ \frac{  M_ \delta (\Theta)}{ S_\delta(\Theta)}  \vEE{\ind{S_\delta(\Theta) <\IF} }  } \lambda_\delta (dr) \\
&=&  \Ex*{ \frac{ M_ \delta (\Theta)}{ S_\delta(\Theta)}   }\\
&=&  \Ex*{ Z(0)\frac{\sup_{t\in  \delta \Z^d} Z(t)}{ S_\delta(Z)}   },
  }
hence the claim follows from \eqref{drejt}.\\
\vEE{  Note that if $\pk{Z(0)>0}=1$,  by \eqref{rinashero} for all $ t\in \R^d$
$$ \Ex{ Z(t)}=\Ex{ Z(t)  \ind{Z(0) >0}}=\Ex{ Z(0) \ind{Z(-t) >0}} =\Ex{Z(0)} $$
and hence  $\pk{Z(t)>0}=\pk{\Theta(t)>0}=1 $ since $\Ex{Z(t)}=\Ex{Z(0)}=1$. Consequently,  applying \eqref{rinashero2} and utilising \eqref{LUK} for all $r\in \R^d, \eta>0$
\bqny{ %\pk{S_0(Z) < \IF, S_\eta(B^r Z) } =
\lefteqn{	\Ex{ Z(-r) \ind{ S_0(Z) < \IF, 0< S_\eta(B^r Z) < \IF}}}\\
     &=& 	\Ex{ \Theta(-r)\ind{ S_0(\Theta) < \IF, 0< S_\eta(B^r\Theta) < \IF}}\\
	&=& 	\Ex{ \ind{ S_0(\Theta) < \IF, 0< S_\eta(\Theta) < \IF}}\\
	&=&	\Ex{\ind{ S_0(\Theta) < \IF}}\\&=& \Ex{ Z(0) \ind{ S_0(Z) < \IF}}\\
	&=&	 \Ex{ Z(-r) \ind{ S_0(Z) < \IF}}\le 1
}
and thus \eqref{laps} follows.}
  The  reason that  the indicator function did  not appear when we applied  \eqref{rinashero2} in the above calculations is that  $\Theta(t)>0$ almost surely for all $t\in \R^d$. %\vEE{EH: Should we omit the last  sentence??}}
\QED 	

\prooftheo{th11} Let $z\in [0,1]$ be fixed. We consider for simplicity the case $d=1$ and set $X_z(t)=\rrd{W}_z(t)-\sigma^2_z(t)/2$.
We show next that for any positive $T$ the function
$$ A_T(z)=\Ex*{ \sup_{t\in [0,T] \cap \delta \Z}  e^{ X_{z}(t)}}$$
 is continuous in $z$ and thus integrable for any $\delta\ge 0$.
First note that by \eqref{norm} % \eqref{eqV}
\bqny{ \Ex*{ \sup_{t\in [0,T] \cap \delta \Z}  e^{ X_{z}(t)}} &=&
1+	\int_{0}^\IF  e^s \pk*{ \sup_{t\in [0,T] \cap \delta \Z} X_{z}(t)> s} ds< \IF.
}
Since $\sigma_z$ determines the covariance function of $X_{z }$, then  \eqref{sz} implies that
fidi's of \kk{$X_{z+h }$} converge weakly to
those of $X_z$  as $h\to 0$. Moreover, for some $\ve>0$ and all $ \abs{h}< \ve$,  by \eqref{kN}
$$ Var\left( X_{z+h} (s)- X_{z+h}(t)\right)= \sigma^2_{z+h}(t-s) \le C \norm{t-s}^\k2{\nu_{0}} $$
for some $C>0$ and all $t,s \in [0,T]$. Consequently, $X_{z+h}(t), t\in [0,T]$ is tight (with respect to $h$) and converges weakly in the  space of
 real-valued continuous functions on $[0,T]$ equipped with the uniform topology, see \cite{Pit20}[Prop 9.7].
Hence   by the continuous mapping theorem for almost all $s \inr$
$$ \lim_{ h \to 0} \pk*{ \sup_{t\in [0,T] \cap \delta \Z} X_{z+h}(t)> s} =
\pk*{ \sup_{t\in [0,T] \cap \delta \Z} X_{z}(t)> s}.$$
Consequently,  for all  $T>0$
\bqny{ \lim_{ h \to 0}\Ex*{ \sup_{t\in [0,T] \cap \delta \Z}  e^{ X_{z+h}(t)}} &=&
	\Ex*{ \sup_{t\in [0,T] \cap \delta \Z}  e^{ X_{z}(t)}} =A_T(z)< \IF
}
implying the claim that $A_T(z), z\in [0,1]$ is  continuous in $z$.  The derivation above rests on the  application of the dominated convergence theorem that can be justified by Borell-TIS inequality. Indeed for all $z\in [0,1]$
$$
\sup_{t\in [0,T] \cap \delta \Z}  e^{ X_{z}(t)}\leq \sup_{t\in[0,T], v\in[0,1]}e^{\rrd{W}_v(t)}
$$
and
\begin{eqnarray*}
\Ex*{ \sup_{t\in [0,T], v\in[0,1]}  e^{ \rrd{W}_{v}(t)}}%&=&\E*{\sup_{t\in [0,S]\, v\in[0,1]}e^ { Y_{v}(t)}}\\
&=&1+\int_0^\infty e^s\pk*{\sup_{t\in[0,T],v\in [0,1]}\rrd{W}_v(t)>s}ds.
\end{eqnarray*}
By Borell-TIS inequality and  \eqref{kN}
$$
\pk*{\sup_{t\in[0,T], v\in [0,1]}\rrd{W}_v(t)>s}\leq e^{-\frac{(s-Q_1)^2}{2Q_2}}
$$
for sufficiently large $s$ and $Q_1,\,Q_2$ positive constants, \kk{which} justifies the limit above.
Since by assumption (\ref{kN}) for \k2{some $C>0$ and all  $z\in [0,1]$
$$ \sigma_z^2(t) \le C (\norm{t}^{\nu_{0}}+ \norm{t}^{\nu_{\infty}}) =:\sigma^2(t),$$
then for the  Gaussian process
$$X(t)= \rrd{W}_0(t)+\rrd{W}_\infty(t)- \sigma^2(t)/2,\ \ t\in [0,T],$$
with \rrd{$W_0,\,W_\infty$} mutually independent centred Gaussian processes with continuous sample paths, stationary increments and variance function
$C\norm{t}^{\nu_{0}},C\norm{t}^{\nu_{\infty}}$, respectively,}
applying \cite{debicki2017approximation}[Thm 3.1]  we obtain for any $z\in [0,1]$
$$ H_{Z_z}^\delta  \le H_{\widetilde{Z}}^\delta \in (0,\IF),$$
\kk{where $\kk{\widetilde{Z}}(t)=\exp(X(t)), \ t\in [0,T]$} and
$$ \Ex*{ \sup_{t\in [0,T] \cap \delta \Z}  e^{ X_{z}(t)}} \le
 \Ex*{ \sup_{t\in [0,T] \cap \delta \Z}  e^{ X(t)}}< \IF.$$
Consequently, by the measurability of $A_T(z), z\in [0,1]$ and the dominated convergence theorem
 $$ \mathcal{H}_{{Z}}^\delta =   \limit{T} T^{-d} \int_0^1 \Ex*{\sup_{ t\in [0,T]^d \cap \delta \Z^d } e^{ X_z (t)}  }    dz =
 \int_0^1 H_{Z_{z}}^\delta dz \le   H^{\delta}_{\widetilde{Z}} < \IF.$$
 \kk{Thus}, in view of  \netheo{th1} we obtain
 $$ \lim_{\delta \downarrow 0} \mathcal{H}_{{Z}}^\delta
 =\int_0^1\lim_{\delta \downarrow 0} H_{Z_z}^\delta dz = \int_0^1 H_{Z_z}^0 dz \ge 0,$$
 \kk{which completes} the proof.  \QED

{ \bf Acknowledgement}:  We are in debt to the reviewers for their comments and suggestions which led to
substantial improvement of the manuscript.  We thank Dmitry Zaporozhets for several important discussions
and the elegant proof of \eqref{ishk}.
Partial support by SNSF Grant 200021-196888 and  NCN Grant No 2018/31/B/ST1/00370 (2019-2022) is kindly acknowledged.
The project is financed by the Ministry of Science and Higher Education in Poland under the programme
"Regional Initiative of Excellence" 2019 - 2022 project number 015/RID/2018/19
total funding amount 10 721 040,00 PLN.

\end{document}